\documentclass[12pt]{article}
\usepackage{amssymb,amsfonts,amsmath, psfrag,eepic,colordvi}
\newcommand{\rmnum}[1]{\romannumeral #1}
\topskip  
\numberwithin{equation}{section}
\newtheorem{theorem}{Theorem}[section]

\newtheorem{conjecture}[theorem]{Conjecture}

\newtheorem{lemma}[theorem]{Lemma}

\begin{document}
\parskip 6pt

\pagenumbering{arabic}
\def\sof{\hfill\rule{2mm}{2mm}}
\def\ls{\leq}
\def\gs{\geq}
\def\SS{\mathcal S}
\def\qq{{\bold q}}
\def\MM{\mathcal M}
\def\TT{\mathcal T}
\def\EE{\mathcal E}
\def\lsp{\mbox{lsp}}
\def\rsp{\mbox{rsp}}
\def\pf{\noindent {\it Proof.} }
\def\mp{\mbox{pyramid}}
\def\mb{\mbox{block}}
\def\mc{\mbox{cross}}
\def\qed{\hfill \rule{4pt}{7pt}}
\def\block{\hfill \rule{5pt}{5pt}}

\begin{center}
{\Large\bf   Ascent Sequences and 3-Nonnesting Set  Partitions  }
\vskip 6mm
\end{center}

\begin{center}
{\small
Sherry H. F. Yan \\[2mm]
 Department of Mathematics, Zhejiang Normal University, Jinhua
321004, P.R. China
\\[2mm]
  huifangyan@hotmail.com
 \\[0pt]
}
\end{center}

\noindent {\bf Abstract.} A sequence $x=x_1 x_2  \ldots x_n $ is said to be an   ascent
sequence of length $n$  if it satisfies $x_1=0$ and $0\leq x_i\leq
asc(x_1x_2\ldots x_{i-1})+1$ for all $2\leq i\leq n$, where $asc(x_1x_2\ldots x_{i-1})$ is the number of ascents in the sequence $x_1x_2\ldots x_{i-1}$. Recently,  Duncan and Steingr\'{\i}msson  proposed the conjecture   that  210-avoiding ascent sequences  of length $n$ are equinumerous with  3-nonnesting set partitions of
$\{1,2,\ldots, n\}$.  In this paper, we  confirm  this conjecture by showing that    210-avoiding ascent sequences  of length $n$  are in bijection with  3-nonnesting   set partitions of
$\{1,2,\ldots, n\}$ via an intermediate structure of growth diagrams for $01$-fillings of Ferrers shapes.

\noindent {\sc Key words}:   ascent sequence, pattern avoiding,  3-nonnesting set partition, growth diagram, $01$-filling of Ferrers shape.

\noindent {\sc AMS Mathematical Subject Classifications}: 05A05,
05C30.


\section{Introduction}
The objective of this paper is to establish a bijection between    210-avoiding ascent sequences  of length $n$  and  3-nonnesting  set partitions of
$\{1,2,\ldots, n\}$.
Let us give an overview of the notation and terminology.

Given a sequence of integers   $x=x_1x_2\ldots x_n$, we say
that the sequence $x$ has an ascent at position $i$ if
$x_i<x_{i+1}$. The number of ascents of $x$ is denoted by $asc(x)$.
A sequence $x=x_1x_2\ldots x_n$ is said to be an {\em ascent
sequence of length $n$} if it satisfies $x_1=0$ and $0\leq x_i\leq
asc(x_1x_2\ldots x_{i-1})+1$ for all $2\leq i\leq n$. Ascent
sequences were introduced by Bousquet-M$\acute{e}$lou et al.
\cite{melon} in their study of $(2+2)$-free posets. Ascent sequences are  closely connected to many other combinatorial structures. Bousquet-M$\acute{e}$lou et al. \cite{melon} constructed bijections
between unlabeled  $(2+2)$-free posets and ascent sequences, between
ascent sequences and permutations avoiding a certain pattern,
between unlabeled $(2+2)$-free posets and a class of involutions
introduced by Stoimenow \cite{Stoi}.  Dukes and Parviainen \cite{Duck} established a bijection between ascent sequences and  upper triangular
matrices with non-negative  integer entries such that all rows and
columns contain at least one non-zero entry.  We call an ascent sequence with no two consecutive equal
entries a {\em primitive} ascent sequence. For example, the sequence $0120234603$ is a primitive ascent sequence of length $10$.

Analogous to   pattern avoidance on permutations, Duncan and Steingr\'{\i}msson \cite{Duncan}
initiated the study of ascent sequences avoiding certain patterns.  For ascent sequences, a {\em pattern } is a word on nonnegative
integers $\{0,1,\ldots, k\}$, where each element appears at least once.  Given an ascent sequence  $x=x_1x_2\ldots x_n$   and a pattern  $\tau=\tau_1\tau_2\ldots \tau_k$, we  say that a subsequence $x_{i_1}x_{i_2}\ldots x_{i_k}$ of $x$ is a pattern of $\tau$ if it is order-isomorphic to $\tau$. If $x$ contains no subsequence of pattern $\tau$, then we say that $x$ is {\em $\tau$-avoiding}.  Denote by $\mathcal{A}_n(\tau)$ and $\mathcal{PA}_n(\tau)$ the set of $\tau$-avoiding  ordinary and primitive ascent sequences, respectively.  Duncan and Steingr\'{\i}msson \cite{Duncan} proved that $|\mathcal{A}_n(\tau)|=C_n$, the $n$th Catalan number, for any $\tau=101, 0101 $ or $021$. Moreover, they proposed  the following  conjecture.
 \begin{conjecture}\label{con}

   210-avoiding ascent sequences  of length $n$ are equinumerous with  3-nonnesting (3-noncrossing) set partitions of
$\{1,2,\ldots, n\}$.
\end{conjecture}

Note that  Mansour and Shattuck \cite{man} recently derived two recurrence relations on the generating function for 210-avoiding ascent sequences. However, as remarked  by  Mansour and Shattuck,  the conjecture is still open.

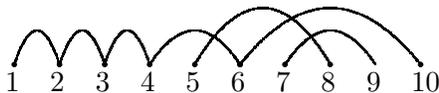
\begin{figure}[h,t]
\begin{center}
\begin{picture}(50,20)
\setlength{\unitlength}{6mm} \linethickness{0.4pt}
 \put(-0.2,0.4){\small$1$}\put(0.8,0.4){\small$2$}
\put(1.8,0.4){\small$3$}\put(2.8,0.4){\small$4$}
 \put(3.8,0.4){\small$5$} \put(4.8,0.4){\small$6$}
\put(5.8,0.4){\small$7$}\put(6.8,0.4){\small$8$}\put(7.8,0.4){\small$9$}\put(8.8,0.4){\small$10$}
\qbezier(0,1)(0.5,2.5)(1,1)\qbezier(1,1)(1.5,2.5)(2,1)\qbezier(2,1)(2.5,2.5)(3,1)\qbezier(3,1)(4,2.5)(5,1)
\qbezier(4,1)(5.5,3.5)(7,1)\qbezier(5,1)(7,3.5)(9,1)\qbezier(6,1)(7,2.5)(8,1)
\put(5,1){\circle*{0.1}} \put(0,1){\circle*{0.1}}\put(1,1){\circle*{0.1}}\put(2,1){\circle*{0.1}}
\put(3,1){\circle*{0.1}}\put(4,1){\circle*{0.1}}\put(5,1){\circle*{0.1}}\put(6,1){\circle*{0.1}}\put(7,1){\circle*{0.1}}
\put(9,1){\circle*{0.1}}

\end{picture}
\vspace{-15pt}
\end{center}
\caption{ The linear representation of  a set partition $\pi=\{\{1,2,3,4, 6,10\}, \{5,8\}, \{7,9\}\}$ .} \label{linear}
\end{figure}

 Recall that a set partition $P$ of $[n] = \{1,2,\ldots, n\}$ can be represented
by a diagram with vertices drawn on a horizontal line in increasing order. For a block
$B$ of $P$, we write the elements of $B$ in increasing order. Suppose that $B=\{i_1, i_2, \ldots, i_k\}$.
Then we draw an arc from $i_1$ to $i_2$, an arc from $i_2$ to $i_3$, and so on. Such a diagram is called
the linear representation of $P$, see Figure \ref{linear} for an example.
 We say that $k$ arcs $(i_1, j_1), (i_2, j_2), \ldots, (i_k, j_k)$ form a
$k$-crossing if $i_1<i_2<\ldots <i_k<j_1<j_2<\ldots<j_k$.  A partition without any $k$-crossing  is said to be {\em k-noncrossing}.
Similarity,  a {$k$-nesting} is a set of $k$ arcs $(i_1, j_1), (i_2, j_2), \ldots, (i_k, j_k)$   such that  $i_1<i_2<\ldots <i_k<j_k<\ldots<j_2<j_1$.
 A set partition without any $k$-nesting  is said to be {\em k-nonnesting}. Chen et al. \cite{chen} proved that k-nonnesting set partitions of $[n]$ are equinumerous with k-noncrossing set
partitions of $[n]$ bijectively using vacillating tableaux as an intermediate object.

In this paper, we aim to  establish a bijection between $210$-avoiding ascents sequences of length $n$ and $3$-nonnesting set partitions of $[n]$
via an intermediate structure of growth diagrams for $01$-fillings of Ferrers shapes described in  \cite{kra} and \cite{van}.

\section{ Growth diagrams for  $01$-fillings of Ferrers shapes}

Our bijection between $210$-avoiding ascent sequences and $3$-nonnesting set partitions
 will be accomplished by the growth diagram for $01$-fillings of Ferrers shapes \cite{kra}. Before we describe
  the growth diagram for $01$-fillings of Ferrers shapes, we give an overview of the notation and terminology.

A {\em partition} $\lambda$ of a positive integer $n$ is defined to
be a sequence $(\lambda_1, \lambda_2, \ldots, \lambda_m)$ of
nonnegative  integers such that $\lambda_1+\lambda_2+\ldots
\lambda_m=n$ and $\lambda_1\geq \lambda_2\ldots \geq \lambda_m$. The empty partition is denoted by $\emptyset$.
For the sake of convenience, we identify a partition $\lambda=(\lambda_1, \lambda_2, \ldots, \lambda_m)$ with the infinite sequence
$(\lambda_1, \lambda_2, \ldots, \lambda_m, 0, 0, \ldots)$, that is, the sequence obtained from $\lambda$ by appending infinitely many $0$'s.
Given a partition $\lambda=(\lambda_1, \lambda_2, \ldots,
\lambda_m)$, its  {\em Ferrers diagram }is the
left-justified array of $\lambda_1+\lambda_2+\ldots+\lambda_m$ squares
with $\lambda_1$ squares in the first row, $\lambda_2$ squares in the
second row, and so on. The {\em partial order} $\subseteq$ on partitions is defined by the containment of their Ferrers diagrams. The {\em conjugate } of a partition $\lambda$ is the partition $(\lambda'_1, \ldots, \lambda'_{\lambda_1})$ where $\lambda'_j$ is the length of the $j$th column in the Ferrers diagram of $\lambda$.

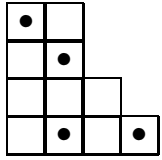
\begin{figure}[h,t]
\begin{center}
\begin{picture}(50,20)
\setlength{\unitlength}{1mm}

\put(0,0){\framebox(5,5) } \put(5,0){\framebox(5,5)  {$\bullet$} }  \put(10,0){\framebox(5,5)  }
 \put(15,0){\framebox(5,5)  {$\bullet$}}

 \put(0,5){\framebox(5,5) } \put(5,5){\framebox(5,5)  }  \put(10,5){\framebox(5,5)  }

 \put(0,10){\framebox(5,5)  } \put(5,10){\framebox(5,5)  {$\bullet$} }

  \put(0,15){\framebox(5,5)  {$\bullet$} } \put(5,15){\framebox(5,5)  }

  \end{picture}
\end{center}
\caption{ Example of a  $01$-filling of a Ferrers shape.  }\label{afilling}
\end{figure}

A Ferrers shape is a Ferrers diagrams in French notation which has straight left side, straight bottom side and supports a descending staircase. We can also encode a Ferrers shape $F$ by sequences of $D's$ and $R's$ by tracing the right/up boundary of $F$ from top-left to bottom-right and writing $D$ (resp. $R$) whenever we encounter a down-step (resp. right-step). For example, the Ferrers shape in Figure \ref{afilling} can be represented  by $RRDDRDRD $.

  A  {\em $01$-filling} of a Ferrers shape is obtained by filling each cells of $F$ with $1's$ and $0's$, see Figure \ref{afilling} for an  example,  where we present $1's$ by $\bullet$ and suppress the $0$'s.
A {\em NE-chain} of a $01$-filling is a sequence of $1's$ such that any $1$ is strictly above and weakly to the right of the preceding $1$ in the sequence.
A {\em SE-chain} of a $01$-filling is a sequence of $1's$ such that any $1$ is weakly below  and strictly to the right of the preceding $1$ in the sequence.
The {\em length} of a $NE$-chain or a $SE$-chain is defined to the number of $1$'s in the chain.

The {\em growth diagram} for a $01$-filling   of Ferrers shape $F$ is obtained by labelling the corners of all the squares in $F$ by partitions in such a way that the partition assigned to any corner  is either equal to the partition to its left or is
obtained from it by adding a {\em horizontal strip}, that is, by a set of squares no two of which are in the same column, and the partition assigned to any corner either equals the
partition below it or is obtained from this partition by adding a {\em vertical  strip}, that is, by a set of squares no two of which are in the same row. We start by assigning the partition $\emptyset$ to each corner on the left and bottom edges of $F$. Then assign the partitions to the other  corners inductively by applying the following forward algorithm.
 Consider the cell in Figure \ref{cell},  filled by $m=0$ or $m=1$ and labeled by the partitions $\rho$, $\mu$, $\upsilon$, where $\rho\subseteq\mu$ and $\rho\subseteq \upsilon$, $\mu$ and $\rho$ differ by a horizontal strip, and $\upsilon$ and $\rho$ differ by a vertical strip. Then $\lambda$ is determined by the following  algorithm.
\begin{itemize}
\item[(F0)] Set $CARRY:=m$ and $i:=1$.
\item[(F1)] Set $\lambda_i=max\{\mu_i+CARRY, \upsilon_i\}$.
\item[(F2)] If $\lambda_i=0$, then stop. The output of the algorithm is $\lambda=(\lambda_1, \lambda_2, \ldots, \lambda_{i-1})$. If not, then set $CARRY:=min\{\mu_i+CARRY, \upsilon_i\}-\rho_i$ and $i:=i+1$. GO to (F1).
    \end{itemize}
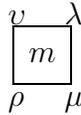
\begin{figure}[h,t]
\begin{center}
\begin{picture}(10,10)
\setlength{\unitlength}{1.5mm}
\put(0,0){\framebox(5,5) {$m$}} \put(-0.5,-2){$\rho$}\put(4.5,-2){$\mu $}\put(-0.5,5.3){$\upsilon$}
\put(4.5,5.3){$\lambda$}

 \end{picture}
\end{center}
\caption{  A cell  filled with $m$.}\label{cell}
\end{figure}

Conversely, given a    labelling of  the corners along the right/up border of a Ferrer shape $F$, one can  reconstruct the labels of the other corners and a $01$-filling of $F$ inductively by applying the following backward algorithm.  Consider the cell in Figure \ref{cell} labelled by the partitions
 $\mu, \upsilon, \lambda$, where $\mu\subseteq\lambda$ and $\upsilon\subseteq \lambda$, where $\lambda$ and $\mu$ differ by a vertical strip, and $\lambda$ and $\upsilon$ differ by a horizontal strip, the backward algorithm works in the following way.
\begin{itemize}
\item[(B0)] Set $i:=max\{j:\lambda_j \,\, \mbox{is positive}\}$ and $CARRY:=0$.
\item[(B1)] Set $\rho_i=min\{\mu_i, \upsilon_i-CARRY\}$.
\item[(B2)] Set $CARRY:=\lambda_i-max\{\mu_i, \upsilon_i-CARRY\}$ and $i:=i-1$. If $i=0$, then stop. The output of the algorithm is $\rho=(\rho_1, \rho_2, \ldots)$ and $m=CARRY$. If not, go to (B1).
 \end{itemize}

  \begin{figure}[h,t]
\begin{center}
\begin{picture}(50,30)
\setlength{\unitlength}{1.5mm}

\put(0,0){\framebox(5,5) } \put(5,0){\framebox(5,5)  {$\bullet$} }  \put(10,0){\framebox(5,5)  }
 \put(15,0){\framebox(5,5)  {$\bullet$}}
    \put(-1.2,-2){$\emptyset$}\put(4.5,-2){$\emptyset$}\put(9.5,-2){$\emptyset$}\put(14.5,-2){$\emptyset$}
\put(19.5,-2){$\emptyset$}
  \put(-1.2,5.3){$\emptyset$}\put(5.3,5.3){$\emptyset$}\put(10.3,5.3){$1$}\put(15.3,5.3){$1$}
\put(20.3,5.3){$11$}

 \put(0,5){\framebox(5,5) } \put(5,5){\framebox(5,5)  }  \put(10,5){\framebox(5,5)  }

    \put(-1.2,10.3){$\emptyset$}\put(5.3,10.3){$\emptyset$}\put(10.3,10.3){$1$}\put(15.3,10.3){$1$}

 \put(0,10){\framebox(5,5)  } \put(5,10){\framebox(5,5)  {$\bullet$} }

  \put(-1.2,15.3){$\emptyset$}\put(5.3,15.3){$\emptyset$}\put(10.3,15.3){$2$}
  \put(0,15){\framebox(5,5)  {$\bullet$} } \put(5,15){\framebox(5,5)  }

  \put(-1.2,20.3){$\emptyset$}\put(5.3,20.3){$1$}\put(10.3,20.3){$21$}

  \end{picture}
\end{center}
\caption{ The growth diagram for the  $01$-filling in Figure \ref{afilling}.}\label{agrowth}
\end{figure}
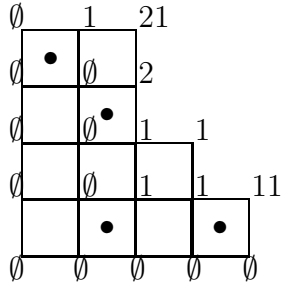

\begin{theorem}{ \upshape   (See {\cite{kra}, Theorem 9})}\label{th1}
Let $F$ be a Ferrers shape given by $D-R$-sequence $w=w_1w_2\ldots w_k$. Then
$01$-fillings of $F$ are in bijection with sequences $(\emptyset=\lambda^0, \lambda^1, \ldots, \lambda^{k}=\emptyset)$ where
$\lambda^{i+1}$ is obtained from $\lambda^{i}$ by doing nothing (i.e., $\lambda^{i+1}=\lambda^{i}$) or  adding  a horizontal strip if $w_i=R$, whereas $ \lambda^{i}$ is obtained from $\lambda^{i-1}$ by doing nothing or deleting a vertical strip if $w_i=D$.
\end{theorem}

\begin{theorem}{ \upshape   (See \cite{kra}, Theorem 10)}\label{th4}
Given a   diagram with empty partitions labelling all the corners along the left side and the bottom side of the Ferrers shape,  suppose that  the corner $c$ is labelled by the partition $\lambda=(\lambda_1, \lambda_2, \ldots, \lambda_l)$. Then $\lambda_1$ is the length of the longest $NE$-chain in the rectangular region to the left and below of $c$, and $\lambda'_1$ is the length of the longest $SE$-chain in the same rectangular region.
\end{theorem}

\begin{lemma}\label{mainle}
Given a $01$-filling of   Ferrers shape $F$ given by $D-R$-sequence $w=w_1w_2\ldots w_k$, let $ (\emptyset=\lambda^0, \lambda^1, \ldots, \lambda^{k}=\emptyset)$ be its corresponding sequence. Then   $\lambda^{i}$ is obtained from $\lambda^{i-1}$ by adding a square   if and only if there is exaclty one $1$ in the column of the cells of $F$ below the corners labeled by $\lambda^{i-1}$ and $\lambda^{i}$, and  $\lambda^{i}$ is obtained from $\lambda^{i-1}$ by deleting  a square   if and only if there is exactlty one $1$ in the row of the cells of $F$ to the left of  the corners labeled by $\lambda^{i-1}$ and $\lambda^{i}$.
\end{lemma}
\pf
Let $c$ be a   cell of $F$ illustrated as Figure \ref{cell}, which is  filled by $m=0$ or $m=1$ and labeled by the partitions $\rho$, $\mu$, $\upsilon$, where $\rho\subseteq\mu$ and $\rho\subseteq \upsilon$, $\mu$ and $\rho$ differ by at most one square, and $\upsilon$ and $\rho$ differ by at most one square.
 In order to prove the lemma, it suffices to show that $\lambda$ has the following properties.
\begin{itemize}
 \item[(a)] If $m=0$ and  $\rho=\upsilon$, then $\lambda=\mu$.
 \item[(b)] If   $m=1$ and $\rho=\upsilon$, then $\lambda$ and  $\mu$ differ by a square.
 \item[(c)] If   $m=0$ and $\rho$ and $\upsilon$ differ by a square, then $\lambda$ and  $\mu$ differ by a square.
 \item[(d)] If $m=0$ and  $\rho=\mu$, then $\lambda=\upsilon$.
 \item[(e)] If   $m=1$ and $\rho=\mu$, then $\lambda$ and  $\upsilon$ differ by a square.
 \item[(f)] If   $m=0$ and $\rho$ and $\mu$ differ by a square, then $\lambda$ and  $\upsilon$ differ by a square.
\end{itemize}

Here  we reformulate the forward algorithm in a slightly different, but of course equivalent fashion. Let $a_1=m$ and $\lambda_1=max\{\mu_1+a_1, \upsilon_1\}$. Suppose that  $\lambda_1, \lambda_2, \ldots, \lambda_i$ and $a_1,a_2, \ldots, a_i$ are already determined. If $\lambda_i=0$, then let $\lambda=(\lambda_1, \lambda_2, \ldots, \lambda_{i-1})$. Otherwise,  let $a_{i+1}=min\{\mu_i+a_i, \upsilon_i\}-\rho_i$ and $\lambda_{i+1}=max\{\mu_{i+1}+a_{i+1}, \upsilon_{i+1}\}$.
Repeat the above procedure until we get $\lambda_k=0$ and let $\lambda=\{\lambda_1, \lambda_2, \ldots, \lambda_{k-1}\}$.

If $m=0$ and $\rho=\upsilon$, then we have $a_i=0$ for all $i\geq 1$ since $\mu_i\geq  \rho_i$ for all $i\geq 1$. This implies that $\lambda=\mu$.

 If $m=1$ and $\rho=\upsilon$, then we have $a_1=1$ and $a_i=0$ for all $i\geq 2$. According to the forward algorithm, we have $\lambda_1=max\{\mu_1+1, \upsilon_1\}=\mu_1+1$ and $\lambda_{i}=max\{\mu_i+a_i, \upsilon_i\}=\mu_i$ for all $i\geq 2$.   This yields that $\lambda$ and  $\mu$ differ by a square.

If $m=0$ and $\rho$ and $\upsilon$ differ by a square, then choose the integer $j$ such that $\upsilon_j=\rho_j+1$. In this case we have $a_i=0$ for all $i\leq  j$ and $\lambda_i=max\{\mu_i+a_i, \upsilon_i\}=\mu_i$ for all $i<j$. Moreover, for all $i\geq j+1$, we have $a_{i+1}=min\{\mu_i+a_i, \upsilon_i\}-\rho_i=0$ and $\lambda_{i+1}=max\{\mu_{i+1}+a_{i+1}, \upsilon_{i+1}\} =\mu_{i+1}$.

 In order to verify property $(c)$, it remains to show that we have either  $\mu_j=\lambda_j$ and $\lambda_{j+1}=\mu_{j+1}+1$, or $\mu_{j+1}=\lambda_{j+1}$ and $\lambda_j=\mu_j+1$. We have two cases.   If $\mu_j>\rho_j$, then we have $\lambda_j=max\{\mu_j+a_j, \upsilon_j\}=max\{\mu_j, \upsilon_j\}=\mu_j$, $a_{j+1}=min\{\mu_j, \upsilon_j\}-\rho_j=1$ and $\lambda_{j+1}=max\{\mu_{j+1}+a_{j+1}, \upsilon_{j+1}\}=\mu_{j+1}+1$ since $\upsilon_{j+1}=\rho_{j+1}\leq \mu_{j+1}$.   If $\mu_j= \rho_j$, then we have $\lambda_j=max\{\mu_j+a_j, \upsilon_j\}=max\{\mu_j, \upsilon_j\}=\upsilon_j=\mu_j+1$,   $a_{j+1}=min\{\mu_j, \upsilon_j\}-\rho_j=0$ and $\lambda_{j+1}=max\{\mu_{j+1}+a_{j+1}, \upsilon_{j+1}\}=\mu_{j+1}$ since $\upsilon_{j+1}=\rho_{j+1}\leq \mu_{j+1}$.  So far,  we have reached the conclusion that $\lambda$ has the properties $(a)-(c)$. By the same reasoning   as in the proofs of $(a)-(c)$, we can verify  properties $(d)-(f)$. The details are omitted. This completes the proof.  \qed

In this paper, we are mainly concerned  with $01$-fillings of   triangular  shape.  Let $\Delta_n$ be the triangular shape with $n$ cells in the bottom row, $n-1$ cells in the row above, etc., and $1$ cell in the topmost
row.
Combining Theorems \ref{th1},  \ref{th4} and Lemma \ref{mainle}, we have the following theorems.

\begin{theorem}\label{th2}
  $01$-fillings of $\Delta_n$ with the property that  every row contains exactly one $1$ are in bijection with sequences $(\emptyset=\lambda^0, \lambda^1, \ldots, \lambda^{2n}=\emptyset)$ where
$\lambda^{2i+1}$ is obtained from $\lambda^{2i}$ by doing nothing or adding  a horizontal strip , whereas $\lambda^{2i+2}$ is obtained from $\lambda^{2i+1}$ by  deleting  a square.
\end{theorem}

\begin{theorem}\label{th3}
  $01$-fillings of $\Delta_n$ with the property that every row and every column  contains at most  one $1$ are in bijection with sequences $(\emptyset=\lambda^0, \lambda^1, \ldots, \lambda^{2n}=\emptyset)$ where
$\lambda^{2i+1}$ is obtained from $\lambda^{2i}$ by doing nothing  or adding  a square, whereas $\lambda^{2i+2}$ is obtained from $\lambda^{2i+1}$ by doing nothing or  deleting  a square.
\end{theorem}

\begin{theorem}\label{th5}
  $01$-fillings of $\Delta_n$  with the property every row contains exactly one $1$ and there is no NE-chain of length $k+1$  are in bijection with sequences $(\emptyset=\lambda^0, \lambda^1, \ldots, \lambda^{2n}=\emptyset)$  where the most number  of columns of any $\lambda^i$ is at most $k$, and
$\lambda^{2i+1}$ is obtained from $\lambda^{2i}$ by doing nothing or  adding a horizontal strip, whereas $\lambda^{2i+2}$ is obtained from $\lambda^{2i+1}$ by   deleting  a square.
\end{theorem}

\begin{theorem}\label{th6}
 $01$-fillings of $\Delta_n$ with the property that every row and every column  contains at most  one $1$  and there is no NE-chain of length $k+1$  are in bijection with sequences $(\emptyset=\lambda^0, \lambda^1, \ldots, \lambda^{2n}=\emptyset)$ where the most number  of columns of any $\lambda^i$ is at most $k$, and
$\lambda^{2i+1}$ is obtained from $\lambda^{2i}$ by doing nothing  or adding  a square, whereas $\lambda^{2i+2}$ is obtained from $\lambda^{2i+1}$ by doing nothing or  deleting  a square.
\end{theorem}

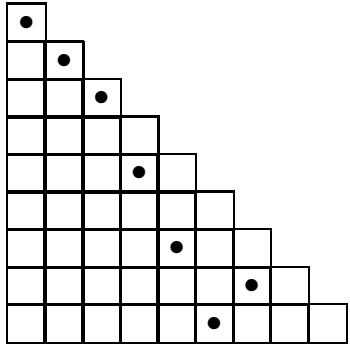
\begin{figure}[h,t]
\begin{center}
\begin{picture}(200,45)
\setlength{\unitlength}{1mm}
\put(60,0){\framebox(5,5) } \put(65,0){\framebox(5,5)  }  \put(70,0){\framebox(5,5)  }
 \put(75,0){\framebox(5,5) }  \put(80,0){\framebox(5,5) }  \put(85,0){\framebox(5,5) {$\bullet$}}
  \put(90,0){\framebox(5,5)  }  \put(95,0){\framebox(5,5)  }  \put(100,0){\framebox(5,5)  }

 \put(60,5){\framebox(5,5)  } \put(65,5){\framebox(5,5)  }  \put(70,5){\framebox(5,5)  }
 \put(75,5){\framebox(5,5) }  \put(80,5){\framebox(5,5) }  \put(85,5){\framebox(5,5)  }
  \put(90,5){\framebox(5,5){$\bullet$} }  \put(95,5){\framebox(5,5)  }

 \put(60,10){\framebox(5,5)  } \put(65,10){\framebox(5,5)  }  \put(70,10){\framebox(5,5)  }
 \put(75,10){\framebox(5,5) }  \put(80,10){\framebox(5,5){$\bullet$}}  \put(85,10){\framebox(5,5)  }
  \put(90,10){\framebox(5,5) }

  \put(60,15){\framebox(5,5)  } \put(65,15){\framebox(5,5)  }  \put(70,15){\framebox(5,5)  }
 \put(75,15){\framebox(5,5)}  \put(80,15){\framebox(5,5) }  \put(85,15){\framebox(5,5)  }

\put(60,20){\framebox(5,5)  } \put(65,20){\framebox(5,5)  }  \put(70,20){\framebox(5,5)  }
 \put(75,20){\framebox(5,5){$\bullet$}}  \put(80,20){\framebox(5,5) }

\put(60,25){\framebox(5,5) } \put(65,25){\framebox(5,5)  }  \put(70,25){\framebox(5,5)  }
 \put(75,25){\framebox(5,5) }

\put(60,30){\framebox(5,5)  } \put(65,30){\framebox(5,5)   }  \put(70,30){\framebox(5,5) {$\bullet$}}

\put(60,35){\framebox(5,5)  } \put(65,35){\framebox(5,5) {$\bullet$}}

\put(60,40){\framebox(5,5) {$\bullet$}}

  \end{picture}
\end{center}
\caption{ A set partition $\pi=\{\{1,2,3,4, 6,10\}, \{5,8\}, \{7,9\}\}$   and its corresponding   $01$-filling.}\label{filling}
\end{figure}

 Recall that there is a bijection between set partitions of $[n]$ and $01$-fillings of $\Delta_{n-1}$ in which every row and every column contains at most one $1$. Given a set partition $\pi$ of $[n]$, we can get a $01$-filling of $\Delta_{n-1}$ by putting a $1$ in the $i$th column and $j$th row from above (where we number rows such that the row consisting of $j-1$ cells is numbered $j$),  whenever $(i,j)$ is an arc in its linear representation. The $01$-filling corresponding to the set partition $\pi=\{\{1,2,3,4, 6,10\}, \{5,8\}, \{7,9\}\}$ is shown in Figure \ref{filling}.    Moreover, a $k$-crossing of a set partition corresponds to a $SE$-chain of length $k$ in the filling, while a $k$-nesting corresponds to a $NE$-chain of length $k$. Thus the following theorem follows immediately.

\begin{theorem}\label{th7}
     $k$-nonnesting set partitions of $[n+1]$  are in bijection with $01$-fillings of $\Delta_n$ with the property that every row and every column  contains at most  one $1$  and there is no NE-chain of length $k$.
\end{theorem}

Denote by $\mathcal{V}_n$ the set of sequences $(\emptyset=\lambda^0, \lambda^1, \ldots, \lambda^{2n}=\emptyset)$  where the most number  of columns of any $\lambda^i$ is at most $2$, and
$\lambda^{2i+1}$ is obtained from $\lambda^{2i}$   by doing nothing or adding a square, whereas $\lambda^{2i+2}$ is obtained from $\lambda^{2i+1}$ by doing nothing or    deleting  a square.

In view of Theorems \ref{th6} and \ref{th7}, in order to prove Conjecture \ref{con}, it suffices to establish a bijection between the set $\mathcal{A}_{n+1}(210)$ and the set $\mathcal{V}_n$. In the next section, we will provide such a bijection.

\section{Proof of the conjecture}
In this section, we will establish a bijection between the set $\mathcal{A}_{n+1}(210)$ and the set $\mathcal{V}_n$.
To this end, we first  construct a bijection between  $210$-avoiding primitive ascent sequences and a certain class of $01$-fillings of triangular shape.

\begin{figure}[h,t]
\begin{center}
\begin{picture}(80,80)
\setlength{\unitlength}{1.5mm}

\put(0,0){\framebox(5,5) } \put(5,0){\framebox(5,5)  }  \put(10,0){\framebox(5,5)  }
 \put(15,0){\framebox(5,5) }  \put(20,0){\framebox(5,5) }  \put(25,0){\framebox(5,5) {$\bullet$}}
  \put(30,0){\framebox(5,5)  }  \put(35,0){\framebox(5,5)  }  \put(40,0){\framebox(5,5)  }
 \put(-1.5,-2){$\emptyset$}\put(4.5,-2){$\emptyset$}\put(9.5,-2){$\emptyset$}\put(14.5,-2){$\emptyset$}
\put(19.5,-2){$\emptyset$}\put(24.5,-2){$\emptyset$}\put(29.5,-2){$\emptyset$}\put(34.5,-2){$\emptyset$}
 \put(39.5,-2){$\emptyset$}\put(44.5,-2){$\emptyset$}
  \put(-1.5,5.3){$\emptyset$}\put(5.3,5.3){$\emptyset$}\put(10.3,5.3){$\emptyset$}\put(15.3,5.3){$\emptyset$}
\put(20.3,5.3){$\emptyset$}\put(25.3,5.3){$\emptyset$}\put(30.3,5.3){$1$}\put(35.3,5.3){$1$}
 \put(40.3,5.3){$1$} \put(45.3,5.3){$1$}
 \put(0,5){\framebox(5,5)  } \put(5,5){\framebox(5,5)  }  \put(10,5){\framebox(5,5)  }
 \put(15,5){\framebox(5,5) }  \put(20,5){\framebox(5,5) }  \put(25,5){\framebox(5,5)  }
  \put(30,5){\framebox(5,5){$\bullet$} }  \put(35,5){\framebox(5,5)  }
  \put(-1.5,10.3){$\emptyset$}\put(5.3,10.3){$\emptyset$}\put(10.3,10.3){$\emptyset$}\put(15.3,10.3){$\emptyset$}
\put(20.3,10.3){$\emptyset$}\put(25.3,10.3){$\emptyset$}\put(30.3,10.3){$1$}\put(35.3,10.3){$2$}
 \put(40.3,10.3){$2$}
 \put(0,10){\framebox(5,5)  } \put(5,10){\framebox(5,5)  }  \put(10,10){\framebox(5,5)  }
 \put(15,10){\framebox(5,5) }  \put(20,10){\framebox(5,5){$\bullet$}}  \put(25,10){\framebox(5,5)  }
  \put(30,10){\framebox(5,5) }
  \put(-1.5,15.3){$\emptyset$}\put(5.3,15.3){$\emptyset$}\put(10.3,15.3){$\emptyset$}\put(15.3,15.3){$\emptyset$}
\put(20.3,15.3){$\emptyset$}\put(25.3,15.3){$1$}\put(30.3,15.3){$11$}\put(35.3,15.3){$21$}
  \put(0,15){\framebox(5,5)  } \put(5,15){\framebox(5,5)  }  \put(10,15){\framebox(5,5)  }
 \put(15,15){\framebox(5,5){$\bullet$}}  \put(20,15){\framebox(5,5) }  \put(25,15){\framebox(5,5)  }
 \put(-1.5,20.3){$\emptyset$}\put(5.3,20.3){$\emptyset$}\put(10.3,20.3){$\emptyset$}\put(15.3,20.3){$\emptyset$}
\put(20.3,20.3){$1$}\put(25.3,20.3){$11$}\put(30.3,20.3){$111$}
\put(0,20){\framebox(5,5)  } \put(5,20){\framebox(5,5)  }  \put(10,20){\framebox(5,5)  }
 \put(15,20){\framebox(5,5){$\bullet$}}  \put(20,20){\framebox(5,5) }
 \put(-1.5,25.3){$\emptyset$}\put(5.3,25.3){$\emptyset$}\put(10.3,25.3){$\emptyset$}\put(15.3,25.3){$\emptyset$}
\put(20.3,25.3){$2$}\put(25.3,25.3){$21$}
\put(0,25){\framebox(5,5) {$\bullet$}} \put(5,25){\framebox(5,5)  }  \put(10,25){\framebox(5,5)  }
 \put(15,25){\framebox(5,5) }
 \put(-1.5,30.3){$\emptyset$}\put(5.3,30.3){$1$}\put(10.3,30.3){$1$}\put(15.3,30.3){$1$}
\put(20.3,30.3){$21$}
\put(0,30){\framebox(5,5)  } \put(5,30){\framebox(5,5)   }  \put(10,30){\framebox(5,5) {$\bullet$}}
\put(-1.5,35.3){$\emptyset$}\put(5.3,35.3){$1$}\put(10.3,35.3){$1$}\put(15.3,35.3){$2$}
\put(-1.5,30.3){$\emptyset$}\put(5.3,30.3){$1$}\put(10.3,30.3){$1$}\put(15.3,30.3){$1$}
\put(20.3,30.3){$21$}
\put(0,35){\framebox(5,5)  } \put(5,35){\framebox(5,5) {$\bullet$}}
\put(-1.5,40.3){$\emptyset$}\put(5.3,40.3){$1$}\put(10.3,40.3){$2$}
\put(0,40){\framebox(5,5) {$\bullet$}}
\put(-1.5,45.3){$\emptyset$}\put(5.3,45.3){$2$}.

  \end{picture}
\end{center}
\caption{ The growth diagram for a $01$-filling of $\Delta_9$.}\label{growth}
\end{figure}
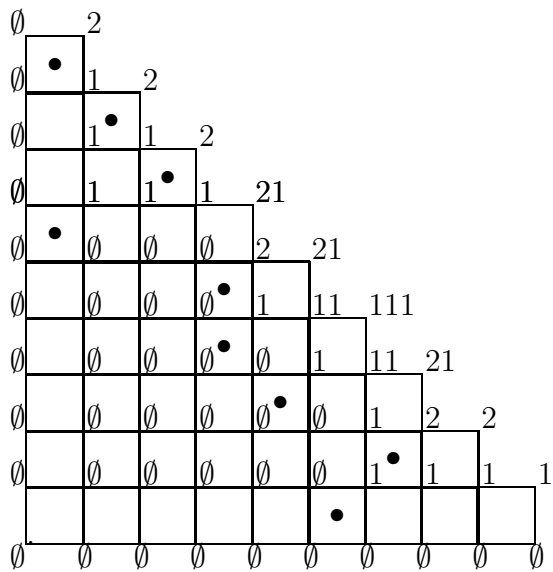

Denote by $\mathcal{N}(\Delta_n)$ the set of all $01$-fillings of $\Delta_n$ with the property that every row contains exactly one $1$ and there is no NE-chain of length $3$. In this section,  a $01$-filling in $\mathcal{N}(\Delta_n)$   will be identified with a sequence $\{( 1, a_1),  (2, a_2),  \ldots, (n, a_n)\}$, where $1\leq a_i\leq i$ and $a_i=k$ if and only if  there is a $1$ in the $i$th row and $k$th column (where we number columns from left to right and  rows such that the row consisting of $j$ cells is numbered $j$).  For example, the filling in figure \ref{growth} is identified with $$\{(1,1), (2,2,), (3,3), (4,1), (5,4), (6,4), (7,5), (8,7), (9,6)\}.$$

Let $n\geq 1$. We aim to describe  a map $\phi$   from the set $\mathcal{PA}_{n+1}(210)$ to the set $\mathcal{N}(\Delta_n)$.     Let $x=x_1x_2\ldots x_{n+1}\in \mathcal{PA}_{n+1}(210)$.  Define $\phi(x)=\{( 1, a_1),  (2, a_2),  \ldots  ,(n, a_n)\}$  where $a_i=i+x_{i+1}-asc(x_1x_2\ldots x_{i+1})$ for all $i=1, 2, \ldots, n$. For example,  let $x=012340415\in \mathcal{PA}_{9}(210)$. Then we have
 $$\phi(x)=\{(1, 1), (2, 2), (3,3), (4,4), (5,1), (6, 5), (7, 3), (8, 7) \}. $$

\begin{lemma}\label{lem1}
  Let  $x=x_1x_2\ldots x_{n+1}\in \mathcal{PA}_{n+1}(210)$ and $\phi(x)=\{( 1, a_1),  (2, a_2),  \ldots  ,(n, a_n)\}$.   Suppose that $a_i\geq a_j$ with $i<j$.  Then we have     $x_{i+1}>x_{j+1}$.
\end{lemma}
\pf    Recall that $a_i=i+x_{i+1}-asc(x_1x_2\ldots x_{i+1})$ and $a_{j}=j+x_{j+1}-asc(x_1x_2\ldots x_{j+1})$.  Since $a_i\geq a_j$, we have
  $$x_{i+1}-x_{j+1}\geq j-i-( asc(x_1x_2\ldots x_{j+1})-asc(x_1x_2\ldots x_{i+1}))\geq 0.$$
  If $x_{i+2}<x_{i+1}$, then we have $asc(x_1x_2\ldots x_{j+1})\leq asc(x_1x_2\ldots x_{i+1})+j-i-1$. This yields that
$$
x_{i+1}-x_{j+1}\geq j-i-( asc(x_1x_2\ldots x_{j+1})-asc(x_1x_2\ldots x_{i+1}))\geq 1.
$$
  Thus we deduce that $x_{i+1}>x_{j+1}$.

  If $x_{i+2}>x_{i+1}$, then we have $asc(x_1x_2\ldots x_{j+1})=asc(x_1x_2\ldots x_{i+1})+1+asc(x_{i+2}\ldots x_{j+1})$. If $asc(x_{i+2}\ldots x_{j+1})=j-i-1$, then we have $x_{i+2}<\ldots <x_{j+1}$.  Thus it follows that $x_{j+1}>x_{i+1}$. This contradicts with the fact that
 $x_{i+1}-x_{j+1}\geq 0$. Hence we have  $ asc(x_{i+2}\ldots x_{j+1})<j-i-1$. It follows that
 $$
 \begin{array}{lll}
 asc(x_1x_2\ldots x_{j+1})&=&asc(x_1x_2\ldots x_{i+1})+1+asc(x_{i+2}\ldots x_{j+1})\\
 &<&asc(x_1x_2\ldots x_{i+1})+j-i.
 \end{array}$$   This implies that $x_{i+1}-x_{j+1}\geq 1$.  This completes the proof. \qed

 \begin{theorem}\label{phi}
 For $n\geq 1$, the map $\phi$ is well defined, that is, for any ascent sequence $x\in \mathcal{PA}_{n+1}(210)$, we have $\phi(x)\in \mathcal{N}(\Delta_n)$.
 \end{theorem}
\pf   Let   $x=x_1x_2\ldots x_{n+1}\in \mathcal{PA}_{n+1}(210)$ and $\phi(x)=\{(1,a_1), (2, a_2), \ldots, (n, a_n)\}$.  First it is necessary to prove that $a_i\leq i$.    According to the definition of ascent sequences, we have $x_{i+1}\leq asc(x_1x_2\ldots x_{i})+1$.

If $x_{i+1}=asc(x_1x_2\ldots x_{i})+1$, then we have $asc(x_1x_2\ldots x_{i+1})=asc(x_1x_2\ldots x_i)+1$. This yields that $x_{i+1}= asc(x_1x_2\ldots x_{i+1})$, which implies that   $a_i=i$.

If $x_{i+1}\leq asc(x_1x_2\ldots x_{i})$, then it follows that $$x_{i+1}\leq asc(x_1x_2\ldots x_i)\leq asc(x_1x_2\ldots x_{i+1}).$$  Recall that $a_i=i+x_{i+1}-asc(x_1x_2\ldots x_{i+1})$. Thus we deduce that  $a_i\leq i$.

Next we aim to show that $a_i\geq 1$. There are  two cases.
If   $asc(x_1x_2\ldots x_{i+1})=i$, then we have $0=x_1<x_2\ldots<x_{i+1}$. In this case, we have
$$
a_i=i+x_{i+1}-asc(x_1x_2\ldots x_{i+1})=  x_{i+1}\geq 1.
$$
If $asc(x_1x_2\ldots x_{i+1})<i$, then we have
$$
a_i=i+x_{i+1}-asc(x_1x_2\ldots x_{i+1})\geq   x_{i+1} +1\geq 1.
$$
Thus we have reached the conclusion that  $1\leq a_i\leq i$. This implies  that
   $\phi(x)$ is  a $01$-filling of $ \Delta_n $ with the property that every row contains exactly one $1$.

It remains to show that  there is no NE-chain   of length $3$ in  $\phi(x)$. Suppose that  the squares $(i, a_i)$ $(j, a_j)$ and $(k,a_k)$ form a NE-chain of length $3$, where $i<j<k$.  This means that $a_i\geq a_j\geq a_k$. By Lemma \ref{lem1}, we have $x_{i+1}>x_{j+1}>x_{k+1}$, which is a $210$ pattern. This leads to a contradiction with the fact that $x$ is $210$-avoiding. Hence, we deduce that there is no NE-chain   of length $3$ in  $\phi(x)$.  This completes the proof. \qed

 In order to show that the map $\phi$ is a bijection,    we proceed   to describe  a map $\phi'$ from the set  $\mathcal{N}(\Delta_n)$ to the set $PA_{n+1}(210)$ for $n\geq 1$.  Given a $01$-filling  $F\in \mathcal{N}(\Delta_n)$, we wish to recover a $210$-avoiding primitive  ascent sequence of length $n+1$. Let $F=\{( 1, a_1),  (2, a_2),  \ldots  ,(n, a_n)\}$. Define $\phi'(F)=(x_1,x_2,\ldots, x_{n+1})$ inductively  as follows:
  \begin{itemize}
  \item   $x_1=0$ and $x_2=1$;
   \item  if $a_{i-1 }<  a_{i }$, then $x_{i+1}=asc(x_1x_2\ldots x_{i})+1+a_{i}-i$  for all $2\leq i\leq n$ ;
 \item if  $a_{i-1}\geq  a_{i}$, then $x_{i+1}=asc(x_1x_2\ldots x_{i})+a_{i}-i$  for all $2\leq i\leq n$.

\end{itemize}
For example, let $F=\{(1, 1), (2, 2), (3,3), (4,4), (5,1), (6, 5), (7, 3), (8, 7)\}\in \mathcal{N}(\Delta_8) $. Then we have $\phi'(F)=012340415$.

\begin{lemma}\label{lem2}
 For any $01$-filling $F=\{( 1, a_1),  (2, a_2),  \ldots  ,(n, a_n)\}\in \mathcal{N}(\Delta_n)$, the sequence $\phi'(F)=x_1 x_2 \ldots x_{n+1}$  has the property that
  $$x_{i+1}=asc(x_1x_2\ldots x_{i+1})+a_{i}-i $$   for all $i\geq 1$. Furthermore,  we have $x_i>x_{i+1}$ if $a_{i-1}\geq a_{i}$, whereas  $x_i<x_{i+1}$ if $a_{i-1}< a_{i}$.
\end{lemma}
\pf   We proceed by induction on $i$. It is easy to check that the statement  holds for $i=1$.   Assume that  $x_{i}=asc(x_1x_2\ldots x_{i})+a_{i-1}-(i-1)$.  Now we shall  show that $x_{i+1}=asc(x_1x_2\ldots x_{i+1})+a_{i}-i$. We consider two cases.

 If $a_{i-1}< a_{i}$, then we have $x_{i+1}=asc(x_1x_2\ldots x_{i})+1+a_{i}-i$. In this case, we have $x_{i+1}> asc(x_1x_2\ldots x_{i})+a_{i-1}-(i-1)=x_i$. This implies   that $x_{i+1}>x_i$ and  $asc(x_1x_2\ldots x_{i+1})=asc(x_1x_2\ldots x_{i})+1$. Thus we deduce that $x_{i+1}=asc(x_1x_2\ldots x_{i+1})+a_{i}-i$ and $x_i<x_{i+1}$.

  If $a_{i-1}\geq  a_{i}$, then we have $x_{i+1}=asc(x_1x_2\ldots x_{i})+a_{i}-i$. It follows that $x_{i}=asc(x_1x_2\ldots x_{i})+a_{i-1}-(i-1)\geq  asc(x_1x_2\ldots x_{i})+1+a_{i}-i=x_{i+1}+1$. This implies that $x_i>x_{i+1}$ and  $asc(x_1x_2\ldots x_{i+1})=asc(x_1x_2\ldots x_{i})$. Thus we have $x_{i+1}=asc(x_1x_2\ldots x_{i+1})+a_{i}-i$ and $x_i>x_{i+1}$.  This completes the proof. \qed

 \begin{theorem}\label{psi}
 Let $n\geq 1$.  The map $\phi'$ is well defined, that is, for any $01$-filling $F\in \mathcal{N}(\Delta_n)$,   we have $\phi'(F)\in \mathcal{PA}_{n+1}(210)$.
 \end{theorem}
 \pf Let $F=\{( 1, a_1),  (2, a_2),  \ldots  ,(n, a_n)\}$ and $\phi'(F)=x_1x_2\ldots x_{n+1}$.
Since $a_i\leq i$, we have $x_{i+1}\leq asc(x_1x_2\ldots x_{i})+1$ according to the definition of the map $\phi'$.
  We next prove  that $x_{i}\geq 0$   by induction on $i$. It is apparent that $x_1\geq 0$ and $x_2\geq 0$. Assume that $x_j\geq 0$ for all $1\leq j\leq i$.  Now we proceed  to show that $x_{i+1}\geq 0$. We have two cases. If $a_{i-1}< a_i$, then we have
  $$
  \begin{array}{ll}
  x_{i+1}&=asc(x_1x_2\ldots x_i)+1+a_i-i\\
     &> asc(x_1x_2\ldots x_i)+a_{i-1}-(i-1).
    \end{array}
  $$
  By Lemma \ref{lem2}, we have $x_{i}=asc(x_1x_2\ldots x_i)+a_{i-1}-(i-1)$. Thus it follows that $x_{i+1}> x_i\geq 0$.

 If $a_{i-1}\geq a_i$, then we have $x_{i+1}=asc(x_1x_2\ldots x_i)+a_i-i$ according to the definition of $\phi'$. If $asc(x_1x_2\ldots x_i)=i-1$, then we have
 $x_{i+1}\geq 0$ since $a_i\geq 1$. Suppose that $asc(x_1x_2\ldots x_i)<i-1$. By Lemma \ref{lem2},    there exists integer $j$ such that $j\leq i-1$ and  $ a_{j-1}\geq a_j$. Let $k$ be the largest such integer.  Hence we have
 $$
 asc(x_1x_2\ldots x_i)=asc(x_1x_2\ldots x_{k+1})+i-k-1.
 $$
Moreover, since  there is no NE-chain of length $3$ in $F$, we have $a_i>a_k$.
 Thus we deduce   that
 $$
  \begin{array}{ll}
  x_{i+1}&=asc(x_1x_2\ldots x_i)+a_i-i\\
     &= asc(x_1x_2\ldots x_{k+1})+a_{i}+(i-k-1)-i\\
     &=asc(x_1x_2\ldots x_{k+1})+a_{i}-k-1\\
     &\geq asc(x_1x_2\ldots x_{k+1})+1+a_{k}-k-1\\
     &=asc(x_1x_2\ldots x_{k+1})+a_{k}-k.
    \end{array}
  $$
 By Lemma \ref{lem2}, we have $x_{k+1}=asc(x_1x_2\ldots x_{k+1})+a_{k}-k$. This yields that $x_{i+1}\geq x_{k+1}\geq 0$.  Hence, we have reached the conclusion that $x_1=0$ and $0\leq x_{i+1}\leq asc(x_1x_2\ldots x_{i})+1$ for all $i\geq 1$.   This ensures that   the obtained sequence $\phi'(F)$ is an ascent sequence of length $n+1$. By Lemma \ref{lem2}, we have $x_{i}\neq x_{i+1}$ for all $i\geq 1$.   Hence the obtained sequence $\phi'(F)$  is primitive.

It remains to show that  $\phi'(F)$ is $210$-avoiding. We claim that for any integers $i,j$ with $i<j$, if $x_{i+1}>x_{j+1}$, then $a_i\geq a_j$. By Lemma \ref{lem2}, we have  $x_{i+1}=asc(x_1x_2\ldots x_{i+1})+a_i-i$ and $x_{j+1}=asc(x_1x_2\ldots x_{j+1})+a_j-j$.
  Suppose that $x_{i+1}>x_{j+1}$. It follows that $asc(x_{i+1}\ldots x_{j+1})\leq j-i-1$.
   Then we have $$
\begin{array}{lll}
a_i-a_j&\geq  & j-i-(asc(x_1x_2\ldots x_{j+1})-asc(x_1x_2\ldots x_{i+1}))-1\\
&=&j-i- asc(x_{i+1}\ldots x_{j+1})-1\\
&\geq& 0.
\end{array}$$  This yields that $a_i\geq a_j$. Hence the claim is proved.

Suppose that  $x_{i+1}>x_{j+1}>x_{k+1}$ with $i<j<k$. Then we have $a_i\geq a_j\geq a_k$. This implies that the squares $(i,a_i)$, $(j,a_j)$ and $(k,a_k)$ form a NE-chain of length $3$, which leads to a contradiction.  Therefore, the obtained sequence $\phi'(F)$ is  $210$-avoiding.  This completes the proof. \qed

\begin{theorem}\label{main}
The maps $\phi$  and $\phi'$  are inverses of each other.
\end{theorem}
\pf
We first prove  that $\phi'(\phi(x))=x$ for any ascent sequence $x=x_1x_2\ldots x_{n+1}\in \mathcal{PA}_{n+1}(210)$.
   Suppose   that $\phi(x)=\{(1,a_1), (2, a_2), \ldots, (n, a_{n})\}$. According to the definition of the map $\phi$, we have $a_i=i+x_{i+1}-asc(x_1x_2\ldots x_{i+1})$. When we apply the map $\phi'$ to $\phi(x)$, we get a 210-avoiding ascent sequence $x'=x'_1 x'_2 \ldots  x'_{n+1}$ such that
\begin{itemize}
  \item   $x'_1=0$ and $x'_2=1$;
   \item  if $a_{i-1 }<  a_{i }$, then $x'_{i+1}=asc(x'_1x'_2\ldots x'_{i})+1+a_{i}-i$ for all $2\leq i\leq n$;
 \item if  $a_{i-1}\geq  a_{i}$, then $x'_{i+1}=asc(x'_1x'_2\ldots x'_{i})+a_{i}-i$ for all $2\leq i\leq n$.
\end{itemize}
We proceed to show that $x'_i=x_i$ for all $i\geq 1$ by induction on $i$. It is clear that $x'_1=0=x_1$ and $x'_2=1=x_2$. Assume that $x'_k=x_k$ for $1\leq  k\leq i$. Now we proceed  to show that $x'_{i+1}=x_{i+1}$. We have two cases.

If $a_{i-1}\geq a_i$, then we have $x'_{i+1}=asc(x'_1x'_2\ldots x'_{i})+a_{i}-i$. By the induction hypothesis, we have $ asc(x'_1x'_2\ldots x'_{i})=asc(x_1x_2\ldots x_{i})$. This implies  that $x'_{i+1}=asc(x_1x_2\ldots x_{i})+a_{i}-i$. By Lemma \ref{lem1}, we have $x_{i}>x_{i+1}$. This yields that $asc(x_1x_2\ldots x_{i+1})=asc(x_1x_2\ldots x_{i})$ and $x_{i+1}=asc(x_1x_2\ldots x_{i+1})+a_{i}-i$. According to the definition of the map $\phi$, we have $a_i=i+x_{i+1}-asc(x_1x_2\ldots x_{i+1})$. Thus we have $x'_{i+1}=x_{i+1}$.

If $a_{i-1}<a_i$, then we have $x'_{i+1}=asc(x'_1x'_2\ldots x'_{i})+1+a_{i}-i$.   By Lemma \ref{lem2}, we have $x'_{i+1}=asc(x'_1x'_2\ldots x'_{i+1})+a_{i}-i$ and  $x'_i<x'_{i+1}$. This implies that $asc(x'_1x'_2\ldots x'_{i+1})=asc(x'_1x'_2\ldots x'_{i})+1$. By the induction hypothesis,  we have  $asc(x'_1x'_2\ldots x'_{i})=asc(x_1x_2\ldots x_{i}) $. It follows that
  $$asc(x'_1x'_2\ldots x'_{i+1})=asc(x'_1x'_2\ldots x'_{i})+1=asc(x_1x_2\ldots x_{i})+1. $$
 By lemma \ref{lem1}, we have $x_{i+1}>x_i$. This yields that $asc(x_1x_2\ldots x_{i+1} )=asc(x_1x_2\ldots x_i)+1$.
    Thus we deduce that $x'_{i+1}=asc(x_1x_2\ldots x_{i+1})+a_{i}-i$. Recall that $a_i=i+x_{i+1}-asc(x_1x_2\ldots x_{i+1})$. Thus we deduce that $x'_{i+1}=x_{i+1}$. Hence we have reached the conclusion that $\phi'(\phi(x))=x$.

  Next we turn to the proof of $\phi(\phi'(F))=F$ for any $01$-filling $F\in \mathcal{N}(\Delta_{n})$. Let $F=\{(1,a_1), (2,a_2), \ldots, (n,a_n)\}$.  Suppose that $\phi'(F)=x_1 x_2 \ldots  x_{n+1}$. When we apply the map $\phi$ to $\phi'(F)$, we get a $01$-filling $F'=\{(1, a'_1), (2, a'_2), \ldots, (n, a'_n)\}\in \mathcal{N}(\Delta_n)$. By the definition of the map $\phi$ and Lemma \ref{lem2}, we have $x_{i+1}=asc(x_1x_2\ldots x_{i+1})+a_i-i$ and $a'_i=i+x_{i+1}-asc(x_1x_2\ldots x_{i+1})$ for $i\geq 1$. This implies that $a'_i=a_i$ for $i\geq 1$. Thus we deduce that $\phi(\phi'(F))=F$. This completes the proof. \qed

Now we are ready to describe a bijection between the  set $\mathcal{A}_{n+1}(210)$ and the set $\mathcal{V}_n$. Our bijection relies on the growth diagrams for $01$-fillings of triangular shape and the bijection $\phi$.

\begin{theorem}\label{mainth}
For $n\geq 0$, there is a bijection between the set $\mathcal{A}_{n+1}(210)$ and the set $\mathcal{V}_n$.
\end{theorem}
\pf   Let $x$ be a $210$-avoiding ascent sequence of length $n+1$, we wish to construct a sequence $\psi(x)\in \mathcal{V}_{n}$. It is apparent that the ascent  sequence
 $x$ can be written   as $x_1^{c_1} x_2^{c_2} \ldots  x_{k+1}^{c_{k+1}}$, where $x_i\neq x_{i+1}$  and $c_i\geq 1$ for  all $i\geq 1$. Let $x'=x_1 x_2\ldots  x_{k+1}$. Obviously, $x'$ is a $210$-avoiding primitive ascent sequence of length $k+1$.
  For $k=0$, let  $\psi(x) $ be a sequence of length $2n+1$ consisting  of  $\emptyset$'s.
  For $k\geq 1$,   by Theorem \ref{phi},   we  obtain        a $01$-filling $\phi(x')$ in $\mathcal{N}(\Delta_{k})$  by applying  the map $\phi$ to $x'$.   When we apply the forward algorithm  to $\phi(x')$, we get the growth diagram for the $01$-filling $\phi(x')$. Let $(\emptyset=\lambda^0, \lambda^1, \ldots, \lambda^{2k}=\emptyset)$ be the sequence of partitions labelling the right/up border of the growth diagram from top-left to bottom-right. By Theorem \ref{th5},  the sequence $(\emptyset=\lambda^0, \lambda^1, \ldots, \lambda^{2k}=\emptyset)$ has the properties that
     the most number  of columns of  any $\lambda^i$ is at most $2$, and
$\lambda^{2i+1}$ is obtained from $\lambda^{2i}$ by doing nothing or adding a  horizontal strip, whereas $\lambda^{2i+2}$ is obtained from $\lambda^{2i}$ by deleting  a square.
 Now we proceed to generate  a sequence of partitions from the sequence $(\emptyset=\lambda^0, \lambda^1, \ldots, \lambda^{2k}=\emptyset)$ by the following procedure.
\begin{itemize}
 \item For all $i\geq 0$, if $\lambda^{2i+1}$ is obtained from $\lambda^{2i}$ by adding two squares, then let $\lambda^{2i+1}=\lambda^{2i+2}$. Let   $(\mu^0, \mu^1, \ldots, \mu^{2k})$ denote  the resulting sequence.
     \item  If $c_1>1$, then adjoin $2(c_1-1)$ copies of $\mu^0$  immediately left to $\mu^{0}$.
\item  For all $1\leq i\leq k$, if $c_{i+1}>1$ ,  then insert $2(c_{i+1}-1)$ copies of  $\mu^{2i}$ immediately left to $\mu^{2i}$ and right to $\mu^{2i-1}$.
  \end{itemize}
Since $c_1+c_2+\ldots+c_{k+1}=n+1$, the resulting sequence is of length $2n+1$.
 Denote by  $( \rho^0, \rho^1, \ldots, \rho^{ 2n})$ the resulting sequence and  set $\psi(x)= ( \rho^0, \rho^1, \ldots, \rho^{ 2n})$.

 In order to show that $\psi(x)\in \mathcal{V}_n$, it suffices to show that  the sequence $( \rho^0, \rho^1, \ldots, \rho^{ 2n})$
 has the following properties:
 \begin{itemize}
 \item[{\upshape (\rmnum{1})}] $\rho^0=\emptyset$ and $\rho^{2n}=\emptyset$;
 \item [{\upshape (\rmnum{2})}]$\rho^{2i+1}$ is obtained from $\rho^{2i}$ by doing nothing  or adding  a square;
 \item [{\upshape (\rmnum{3})}]$\rho^{2i+2}$ is obtained from $\rho^{2i+1}$ by doing nothing or deleting a square;
 \item [{\upshape (\rmnum{4})}]each $\rho^i$ has at most two columns.
    \end{itemize}
 From the construction of the map $\psi$, it is easily seen that the sequence $( \rho^0, \rho^1, \ldots, \rho^{ 2n})$ has properties {\upshape (\rmnum{1})} and {\upshape (\rmnum{4})}. Moreover,  it is easy to check  that the sequence $(\rho^0, \rho^1, \ldots, \rho^{2n})$ is defined by
\begin{equation}\label{eq1}
\rho^l=\left \{ \begin{array}{ll}
\mu^{2i},& \,\, \mbox{if} \,\, \sum_{j=0}^{i}2c_j\leq l\leq \sum_{j=0}^{i+1}2c_{j}-2,\\
\mu^{2i+1},& \,\, \mbox{if} \,\, l=\sum_{j=0}^{i+1}2c_{j}-1,
\end{array}\right.
\end{equation}
with the assumption $c_0=0$.

 We claim that we have either $(1)$ $\mu^{2i+1}$ is obtained from $\mu^{2i}$ by doing nothing or adding a square and $\mu^{2i+2}$ is obtained from $\mu^{2i+1}$ by deleting a square, or $(2)$ $\mu^{2i+1}$ is obtained from $\mu^{2i}$ by  adding a square and $\mu^{2i+2}$ is obtained from $\mu^{2i+1}$ by doing nothing.      Since each $\lambda^i$ has at most two columns,  $\lambda^{2i+1}$ is obtained from $\lambda^{2i}$ by adding   a horizontal strip of at most two squares.
 If  $\lambda^{2i+1}$ is obtained from $\lambda^{2i}$ by  doing nothing or adding a square, then we have $\mu^{2i}=\lambda^{2i}, \mu^{2i+1}=\lambda^{2i+1}$, and $\mu^{2i+2}=\lambda^{2i+2}$. Recall that $\lambda^{2i+2}$ is obtained from $\lambda^{2i+1}$ by deleting a square.  Thus  $\mu^{2i+1}$ is obtained from $\mu^{2i}$ by doing nothing or   adding a square and $\mu^{2i+2}$ is obtained from $\mu^{2i+1}$ by  deleting a square.

If  $\lambda^{2i+1}$ is obtained from $\lambda^{2i}$ by adding a horizontal strip of two squares, then according to the construction of the map $\psi$,  we have  $\mu^{2i+1}=\lambda^{2i+2}=\mu^{2i+2}$.     Recall that $\lambda^{2i+2}$ is obtained from $\lambda^{2i+1}$ by deleting a square. Hence $\lambda^{2i+2}$ and  $\lambda^{2i}$  differ by exactly one square. Since $\mu^{2i}=\lambda^{2i}$ and $\mu^{2i+2}=\lambda^{2i+2}=\mu^{2i+1}$,    it follows that $\mu^{2i+1}$ is obtained from $\mu^{2i}$ by adding a square and $\mu^{2i+2}$ is obtained from $\mu^{2i+1}$ by doing nothing.       Hence the claim is proved.

      By Formula \ref{eq1}, we deduce     that either $(1)$ $\rho^{\sum_{j=0}^{i+1}2c_j-1}$ is obtained form $\rho^{\sum_{j=0}^{i+1}2c_j-2}$ by doing nothing or adding a square and $\rho^{\sum_{j=0}^{i+1}2c_j}$ is obtained form $\rho^{\sum_{j=0}^{i+1}2c_j-1}$ by  deleting  a square, or $(2)$ $\rho^{\sum_{j=0}^{i+1}2c_j-1}$ is obtained form $\rho^{\sum_{j=0}^{i+1}2c_j-2}$ by   adding a square and $\rho^{\sum_{j=0}^{i+1}2c_j}$ is obtained form $\rho^{\sum_{j=0}^{i+1}2c_j-1}$ by doing nothing.
      Moreover, for all $\sum_{j=0}^{i}2c_j< l\leq \sum_{j=0}^{i+1}2c_j-2$, we have  $\rho^{l}=\rho^{l-1}$.  Hence,  the sequence $( \rho^0, \rho^1, \ldots, \rho^{ 2n})$ has  properties {\upshape (\rmnum{2})} and {\upshape (\rmnum{3})}, which implies  that $\psi(x)\in \mathcal{V}_n$.

  Conversely, given a  sequence $V=( \rho^0, \rho^1, \ldots, \rho^{ 2n})\in \mathcal{V}_n$, we wish to recover a $210$-avoiding ascent sequence $\psi'(V)$.
   If $\rho^i=\emptyset$ for all $i\geq 0$, then let $\psi'(V)$ be the ascent sequence of length $n+1$ consisting of $0$'s.
   Otherwise,    for all $i\geq 0$, remove  $\rho^{2i}$ and $\rho^{2i+1}$ from the sequence $( \rho^0, \rho^1, \ldots, \rho^{ 2n})$  whenever  $\rho^{2i}=\rho^{2i+1}=\rho^{2i+2}$.
  Assume that
 $$(\rho^{2j_0}, \rho^{2j_0+1}, \rho^{2j_1}, \rho^{2j_1+1}, \ldots, \rho^{2j_{m-1}}, \rho^{2j_{m-1}+1}, \rho^{2j_{m}})$$ is the resulting sequence.
 Let $(\upsilon^0, \upsilon^1, \ldots, \upsilon^{2m})$ be a sequence such that $(1)$ $\upsilon^{2i}=\rho^{2j_i}$; $(2)$
 if   $\rho^{2j_i+1}=\rho^{2j_{i+1}}$, then $\upsilon^{2i+1}$ is a partition obtained from
 $\upsilon^{2i}$ by adding a horizontal strip of two squares; $(3)$ otherwise, $\upsilon^{2i+1}=\rho^{2j_i+1}$.
 By  applying  the backward algorithm to $(\upsilon^0, \upsilon^1, \ldots, \upsilon^{2m})$, we obtain a $01$-filling $F\in \mathcal{N}(\Delta_{m})$.  Let $\phi'(F)=y_1 y_2 \ldots  y_{m+1}$ and $\psi'(V)=y^{c'_1}_1 y^{c'_2}_2 \ldots  y_{m+1}^{c'_{m+1}}$, where $c'_{i+1}=j_{i}-j_{i-1}$ for $i\geq 1$  and $c'_1=j_0+1$.

 It is apparent that we have $j_{m}=n$.
Since $c'_1+c'_2+\ldots +c'_{m+1}=n+1$, the obtained sequence $\psi'(V)$ is a sequence of length $n+1$.
By Theorems \ref{th5} and \ref{psi}, in order to show that $\psi'(V)\in \mathcal{A}_{n+1}(210)$, it
suffices to show that  the sequence $(\upsilon^0, \upsilon^1, \ldots, \upsilon^{2m})$ verifies the following points.
\begin{itemize}
 \item[{\upshape (\rmnum{1})}$^{'}$] $\upsilon^0=\upsilon^{2m}=\emptyset$;
  \item [{\upshape (\rmnum{2})}$^{'}$] $\upsilon^{2i+1}$ is obtained from $\upsilon^{2i}$ by doing nothing or adding a horizontal strip;
   \item [{\upshape (\rmnum{3})}$^{'}$] $\upsilon^{2i+2}$ is obtained from $\upsilon^{2i+1}$ by deleting a square;
   \item [{\upshape (\rmnum{4})}$^{'}$]  each $\upsilon^i$ has at most two columns.
\end{itemize}

It is easy to check that  the statements of   {\upshape (\rmnum{1})}$^{'}$,  {\upshape (\rmnum{2})}$^{'}$ and  {\upshape (\rmnum{4})}$^{'}$ are true  for the sequence $(\upsilon^0, \upsilon^1, \ldots, \upsilon^{2m})$.   From the construction of the map $\psi'$, it is easily seen that if $\rho^{2j_i+1}\neq \rho^{2j_{i+1}}$, then  $\upsilon^{2i+2}=\rho^{2j_{i+1}}$ and
$\upsilon^{2i+1}=\rho^{2j_{i}+1}$. Moreover, we have $\rho^{2j_i+2}=\rho^{2j_{i+1}}$.  This implies that $\rho^{2j_i+2}$ is obtained form $\rho^{2j_i+1}$ by deleting a square. Hence, it follows that   $\upsilon^{2i+2}$ is obtained from $\upsilon^{2i+1}$ by deleting a square.

 In order to verify {\upshape (\rmnum{3})}$^{'}$, it remains to  show that if $\rho^{2j_i+1}=\rho^{2j_{i+1}}$, then $\upsilon^{2i+2}$ is obtained from $\upsilon^{2i+1}$ by deleting a square. Suppose that $\rho^{2j_i+1}=\rho^{2j_{i+1}}$.
 From the construction of the sequence $$(\rho^{2j_0}, \rho^{2j_0+1}, \rho^{2j_1}, \rho^{2j_1+1}, \ldots, \rho^{2j_{m-1}}, \rho^{2j_{m-1}+1}, \rho^{2j_{m}}),$$ it follows that  $\rho^{2j_i+1}\neq \rho^{2j_{i}}$. This implies that $\rho^{2j_i+1}$ is obtained from $\rho^{2j_{i}}$ by adding one square. Recall that $\upsilon^{2i}=\rho^{2j_i}$. This implies that  $\upsilon^{2i+1}$ is obtained from $\upsilon^{2i}$ by adding a horizontal strip of two squares. Since $\upsilon^{2i+2}=\rho^{2j_{i+1}}=\rho^{2j_i+1}$ and $\rho^{2j_i+1}$ is obtained from $\rho^{2j_{i}}$ by adding one square,    the partition  $\upsilon^{2i+2}$ is obtained from $\upsilon^{2i}$ by adding one square.  Recall that $\upsilon^{2i+1}$ is obtained from $\upsilon^{2i}$ by adding a horizontal strip of two squares.  It follows that   $\upsilon^{2i+2}$ is obtained from $\upsilon^{2i+1}$ by deleting a square.
Hence, we have reached the conclusion that $\psi'(V)\in \mathcal{A}_{n+1}(210)$.

Now we proceed to show that the map $\psi$ is indeed a bijection. To this end, we will show that the maps $\psi$ and $\psi'$ are inverses of each other.
We first show that   the map   $\psi'$ is the inverse of the map $\psi$,  that is, $\psi'(\psi(x))=x$.
  By Theorems \ref{main} and \ref{th5}, it suffices to prove that
\begin{equation}\label{eq2}
(\rho^{2j_0}, \rho^{2j_0+1}, \rho^{2j_1}, \rho^{2j_1+1}, \ldots, \rho^{2j_{m-1}}, \rho^{2j_{m-1}+1}, \rho^{2j_{m}})=(\mu^0, \mu^1, \ldots, \mu^{2k}),
\end{equation}
and $c_j=c'_j$ for all $j\geq 1$.
Recall that we have    either $(1)$ $\rho^{\sum_{j=0}^{i+1}2c_j-1}$ is obtained form $\rho^{\sum_{j=0}^{i+1}2c_j-2}$ by doing nothing or adding a square and $\rho^{\sum_{j=0}^{i+1}2c_j}$ is obtained form $\rho^{\sum_{j=0}^{i+1}2c_j-1}$ by  deleting  a square, or $(2)$ $\rho^{\sum_{j=0}^{i+1}2c_j-1}$ is obtained form $\rho^{\sum_{j=0}^{i+1}2c_j-2}$ by   adding a square and $\rho^{\sum_{j=0}^{i+1}2c_j}$ is obtained form $\rho^{\sum_{j=0}^{i+1}2c_j-1}$ by doing nothing.
      Moreover, for all $\sum_{j=0}^{i}2c_j< l\leq \sum_{j=0}^{i+1}2c_j-2$, we have  $\rho^{l}=\rho^{l-1}$. From the construction of the map $\psi'$,
  it follows that $j_i=\sum_{j=0}^{i+1}c_j-1$  for all $i\geq 0$. This implies that
  $c'_1=j_0+1=c_1$ and $c'_{i+1}=j_i-j_{i-1}=c_{i+1}$ for all $i\geq 1$. Moreover,
    by Formula \ref{eq1}, we have $\rho^{\sum_{j=0}^{i+1}2c_j-2}=\mu^{2i}$ and $\rho^{\sum_{j=0}^{i+1}2c_j-1}=\mu^{2i+1}$.  Thus we obtain Formula \ref{eq2}.

     By Theorems \ref{main} and \ref{th5}, it is routine to check that the map $\psi$  reverses each step of the map $\psi'$. This implies that the map $\psi$ is the inverse map of $\psi'$.  
    Hence,   the maps $\psi$ and $\psi'$ are inverses of each other.  
         This completes the proof. \qed

For example, let $x=0012303222353$ be   a $210$-avoiding ascent sequence of length $13$. Then $x$ can be written as $0^2 1^12^13^10^13^12^33^15^13^1$. Let $x'=0123032353$, which is a $210$-avoiding primitive ascent sequence of length $10$.
By applying the map $\phi$ to $x'$, we get a $01$-filling
$$\phi(x')=\{(1,1), (2,2,), (3,3), (4,1), (5,4), (6,4), (7,5), (8,7), (9,6)\}\in \mathcal{N}(\Delta_9).$$  Then we get the growth diagram for the $01$-filling $\psi(x')$ illustrated in Figure \ref{growth}.  Let
$$ (\lambda^0, \lambda^1, \ldots, \lambda^{18} )=(\emptyset, 2, 1,2,1,2,1,21,2,21,11,111,11,21,2,2,1,1,\emptyset) $$ be the sequence of partitions labelling the corners  along the right/up border of the growth diagram from top-left to bottom-right. It is easy to check that
$\lambda^1$ is obtained from $\lambda^0$ by adding two squares and   $\lambda^7$ is obtained form $\lambda^6$ by adding two squares. Hence we get a sequence
$$
(\mu^0, \mu^1, \ldots, \mu^{18} )=(\emptyset, 1, 1,2,1,2,1,2,2,21,11,111,11,21,2,2,1,1,\emptyset)
$$
by replacing $\lambda^1$  and $\lambda^7$ with the partitions $\lambda^2=1$ and $\lambda^8=2$, respectively.
 Finally, we obtain   a sequence  $$\psi(x)=(\emptyset, \emptyset, \emptyset, 1, 1,2,1,2,1,2,2,21,11,  111,11,11,11,11,11,21,2,2,1,1,\emptyset)\in \mathcal{V}_{12} $$
 from  $(\mu^0, \mu^1, \ldots, \mu^{18} )$
  by inserting   two copies  of $\mu^0$ immediately left to $\mu^0$, and inserting four copies of $\mu^{12}$ immediately left to $\mu^{12}$ and right to $\mu^{11}$.

Conversely,  given   a sequence of partitions
$$
V=(\rho^0, \rho^1, \ldots, \rho^{24})=(\emptyset, \emptyset, \emptyset, 1, 1,2,1,2,1,2,2,21,11,  111,11,11,11,11,11,21,2,2,1,1,\emptyset),
$$  we aim to recover an ascent sequence $\psi'(V)$.
It is easily seen that $\rho^0=\rho^1=\rho^2$, $\rho^{14}=\rho^{15}=\rho^{16}$ and $\rho^{16}=\rho^{17}=\rho^{18}$. We obtain a sequence
$$
(\rho^2, \rho^3, \rho^4, \rho^5, \rho^6, \rho^7, \rho^8, \rho^9, \rho^{10}, \rho^{11}, \rho^{12}, \rho^{13}, \rho^{18}, \rho^{19}, \rho^{20}, \rho^{21}, \rho^{22}, \rho^{23}, \rho^{24})
$$
by removing $\rho^0, \rho^1, \rho^{14}, \rho^{15}, \rho^{16}, \rho^{17}$.
Moreover,  we have $\rho^3=\rho^4$ and $\rho^9=\rho^{10}$.
Replace  $\rho^3$  and $\rho^9$ with the partitions $2$ and  $21$ respectively.   This leads to a sequence of partitions
$$ (\upsilon^0, \upsilon^1, \ldots, \upsilon^{18} )=(\emptyset, 2, 1,2,1,2,1,21,2,21,11,111,11,21,2,2,1,1,\emptyset).$$
Applying the backward algorithm, we get a $01$-filling $$F=\{(1,1), (2,2,), (3,3), (4,1), (5,4), (6,4), (7,5), (8,7), (9,6)\}.$$
Finally we have $\phi'(F)=0123032353$ and $\psi'(V)=0012303222353$.

Combining Theorems  \ref{th6}, \ref{th7} and \ref{mainth}, we get the following theorem, which leads to a combinatorial proof of Conjecture \ref{con}.
\begin{theorem}
$210$-avoiding ascent sequences of length $n$ are in bijection with $3$-nonnesting set  partitions of $[n]$.
\end{theorem}

\noindent{\bf Acknowledgments.}   The   author was supported by the
National Natural Science Foundation of China (No.10901141).


\end{document}